\documentclass[pdflatex,sn-mathphys-num]{sn-jnl}

\usepackage[LGR,T1]{fontenc}   
\usepackage[greek.polutoniko,english]{babel}
\usepackage{graphicx}
\usepackage{float} 
\usepackage{amsmath,amssymb,amsfonts}
\usepackage{textcomp}
\usepackage{xurl}

\let\orcidlogo\relax
\usepackage{orcidlink}

\hypersetup{
  pdftitle={Who Wrote This? Turing, Total Variation, and the Mathematics of AI-Text Detection},
  pdfauthor={Santiago Schnell}
}

\providecommand{\burl}[1]{}
\renewcommand{\burl}[1]{\url{#1}}

\begin{document}

\title[Who Wrote This?]{Who Wrote This? Turing, Total Variation, and the Mathematics of AI-Text Detection}

\author*[1,2,3]{%
  \fnm{Santiago} \sur{Schnell}\,
  \orcidlink{0000-0002-9477-3914}%
}
\email{schnells@umich.edu}

\affil*[1]{\orgdiv{Department of Molecular \& Integrative Physiology}, \orgname{University of Michigan Medical School},
\orgaddress{\city{Ann Arbor}, \state{Michigan}, \postcode{48109-5622}, \country{USA}}}

\maketitle

\begin{flushright}
\begin{minipage}{0.78\textwidth}
\raggedleft
\foreignlanguage{greek}{ἐὰν δέ τι ἔρῃ τῶν λεγομένων βουλόμενος μαθεῖν, ἕν τι σημαίνει μόνον ταὐτὸν ἀεί.}\\[3pt]
\emph{\ldots but if you question anything that has been said because you want to learn more, it continues to signify just 
that very same thing forever.}\\[3pt]
--- Plato, \emph{Phaedrus} 275d~\cite{Plato1997}
\end{minipage}
\end{flushright}
\vspace{1.5ex}

\section*{A suspiciously exact number}
Imagine the following scenario. An advanced undergraduate in a biology class submits an essay as a homework assignment. The 
arguments are logical, the prose is fluent, the paragraphs orderly, perhaps a little too orderly. An instructor pastes the text into an 
artificial intelligence (AI) detector. A few seconds later, a number appears:
\[
98\%\ \text{AI-generated}.
\]

The number looks scientific. For a classroom assignment, though, that appearance can become an accusation remarkably quickly if the 
essay was supposed to be written without AI assistance. But what does the 98\% actually measure? Is it meant to represent the 
probability that the essay was generated by AI\@? The probability that the student used AI at some point in producing it? The fraction 
of the text classified as AI-like? Or merely a model score transformed until it resembles a probability? More basically, what can the final
 arrangement of words tell us about the process that produced them?

A 2023 experiment illustrates the difficulty. Weixin Liang and colleagues~\cite{Liang2023} tested seven widely used AI-text detectors on 
91 English-language TOEFL essays written by non-native speakers. Every essay in the collection was human-written. Even so, the 
detectors produced an average false-positive rate of 61.3\%. Eighteen essays were labeled AI-generated by every detector, and 89 were
 flagged by at least one.

The researchers then asked ChatGPT to enrich the vocabulary of those same essays. After this intervention, the reported average 
false-positive rate fell from 61.3\% to 11.6\%. The irony is striking: introducing actual AI assistance made the essays appear more human 
to the detectors.

These findings concern the models and detectors tested in 2023 and do not necessarily generalize to current systems. But the underlying 
inferential problem cannot be dismissed so easily. Given only a finished text of finite length, can we determine, with rigor, reproducibility, 
and robustness, whether it was written entirely by a human, generated by AI, or produced through some mixture of the two?

A second problem lies behind the first. Classifying a text as more human-like or machine-like, attributing it to the process that produced 
it, and certifying that a particular person wrote it under particular conditions are related tasks, but they are not equivalent. They require
different kinds of evidence, and confusing them can turn statistical resemblance into an accusation about a person. Alan Turing gives us 
a useful place to begin.

\section*{Turing's game, played backward}
In 1950, Turing opened his paper ``Computing Machinery and Intelligence''~\cite{Turing1950} with the question ``Can machines think?''
He did not try to answer it directly. The meanings of \emph{machine} and \emph{think}, he argued, were too uncertain to support the inquiry in that form. 
Instead, he replaced the question with one stated in terms of observable behavior: the imitation game.

Turing's original game involved a man, a woman, and an interrogator communicating at a distance through written messages. He then asked what would 
happen if a machine took the place of one of the participants. Could the interrogator still make the correct identification? In the formulation that later 
became known as the Turing test, the essential idea is familiar: an interrogator exchanges messages with unseen respondents and tries to distinguish a 
machine from a human.

The move is important. Turing set aside a disputed internal property---thinking---and replaced it with an experiment based on observable behavior. The 
interrogator need not decide what thinking \emph{is}; she has to decide whether the evidence available to her permits a reliable distinction. Bradford and 
Wollowski~\cite{Bradford1995} later pushed this idea into explicitly mathematical territory by formalizing the Turing test as an interactive proof system and 
studying its limitations using ideas from probability and computational complexity.

AI-text detection poses a related problem, but with a crucial change. There is no interlocutor. There is only a completed object: an essay, a report, a proof, 
this paper, perhaps a poem. From that fixed sequence of words, a detector tries to infer something about an unobserved source:
\[
\text{human or machine?}
\]

Shao, Uchendu, and Lee~\cite{Shao2019} use the phrase ``reverse Turing test'' for precisely the task of using a machine classifier to distinguish human-made 
from machine-made texts. But for our purposes there is a second reversal, and it matters mathematically. Turing's interrogator participates in an 
\emph{interactive} experiment. She can ask a question, observe the answer, choose another question in response, return to an earlier point, and probe 
whatever seems most informative. The evidence can be gathered adaptively.

A detector examining a submitted essay cannot do this. The experiment, if we can call it one, is already over. The words have been chosen, their order is 
fixed, and no new observation can be requested from the process that produced them.

The detector therefore faces a problem of reconstruction rather than interrogation. It observes an outcome and tries to infer a hidden history from the 
traces preserved in that outcome. This is an inverse problem: from a finite string of words, it attempts to reason backward toward the process that 
generated them. How much of that history can a finite string actually preserve?

\section*{The overlap no detector can erase}
To set the stage, we boil the detection question down to its simplest form. The scenario will be unrealistic, but deliberately so: simplifying a problem 
is often the first step toward understanding what makes it difficult.

For now, we set aside mixed human–AI provenance and suppose that every text comes from exactly one of two sources. Let \(\mathcal X\) be the collection 
of all possible finite essays, and imagine two lotteries on that set. The human lottery assigns probability \(P_H(x)\) to each possible text \(x\); the AI lottery 
assigns probability \(P_A(x)\). These are thought-experiment probabilities, not distributions that we could hope to estimate exactly in practice. We observe 
the essay, not the lottery that produced it, and must infer which source it came from.

Now give the detector every possible advantage. Suppose that the two sources are equally likely in advance and that both distributions are known perfectly. 
After observing \(x\), the optimal rule is immediate: choose human if \(P_H(x)>P_A(x)\), and AI if the inequality is reversed; ties may be broken arbitrarily.

Even this ideal detector makes mistakes whenever the two source distributions overlap. At a given text \(x\), the detector chooses the more 
likely source, so the error comes from the less likely one. With equal prior probabilities, the contribution 
of \(x\) to the average error is
\[
\frac12\min\{P_H(x),P_A(x)\}.
\]
Summing over all possible texts gives the smallest achievable error,
\[
\frac12\sum_{x\in\mathcal X}
\min\{P_H(x),P_A(x)\}.
\]

The quantity inside the factor of one-half is the overlap between the two distributions. Total variation measures the complementary separation. Because
\[
\min\{a,b\}=\frac12(a+b-|a-b|)
\]
and each distribution sums to one,
\[
\operatorname{TV}(P_H,P_A)
=\frac12\sum_{x\in\mathcal X}|P_H(x)-P_A(x)|
=1-\sum_{x\in\mathcal X}\min\{P_H(x),P_A(x)\}.
\]

Hence, with equal prior probabilities, the smallest average error of any detector is
\begin{equation}
R_{\mathrm{best}}
=\frac12\left[1-\operatorname{TV}(P_H,P_A)\right].
\label{eq:bayes-risk}
\end{equation}
This identity is standard in statistical decision theory \cite{LeCam1986}. Here it gives us the central constraint: the overlap belongs to the statistical experiment, 
not to the ingenuity of the classifier.

If total variation is zero, the two distributions are identical and no detector can do better than a coin toss. If it lies strictly between zero and one, classification 
may be useful, perhaps extremely useful, but some error remains unavoidable. Zero error requires
\[
\operatorname{TV}(P_H,P_A)=1.
\]
In the present discrete setting, this means that the two distributions have disjoint supports. In probabilistic terms, they 
are mutually singular. In plain language, perfect detection therefore requires separate worlds of writing: every text that can occur under one source must be 
impossible under the other.

If some text \(x\) has positive probability under both processes, then observing \(x\) cannot reveal its source with certainty. One might expect conventional 
phrases, standard transitions, or formulaic prose to contribute heavily to such overlap, but the mathematics itself does not tell us which texts occupy it. That
is an empirical question. The identity says something more basic: wherever overlap exists, no cleverer classifier can recover information that the observation 
does not contain.

\begin{figure}[H]
\centering
\includegraphics[width=0.94\textwidth]{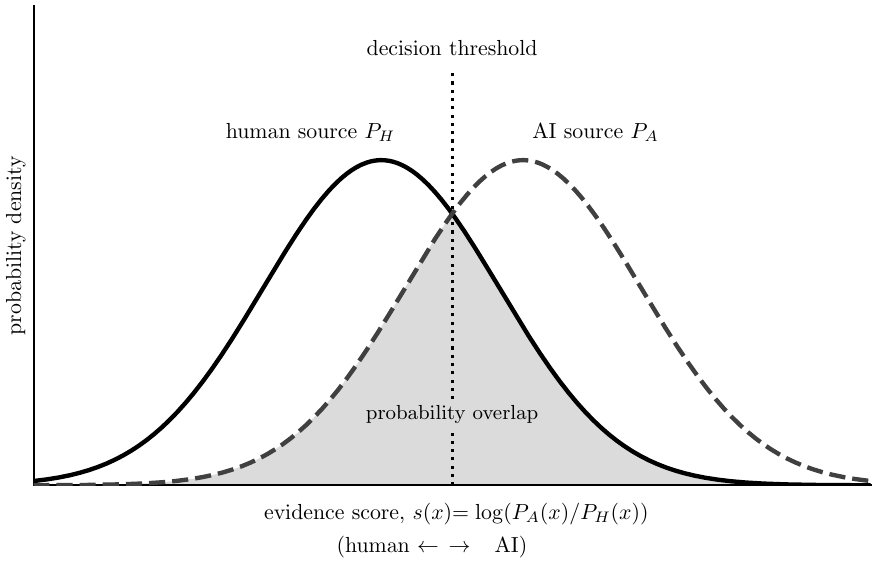}
\caption{\textbf{The overlap no detector can erase.} The horizontal axis is the evidence score \(s(x)\), increasing from evidence favoring the human source to evidence favoring the AI source; the curves show its probability densities under the two sources in a Gaussian example. The shaded region is the overlap between these densities. Moving the decision threshold trades false positives against false negatives but cannot eliminate this overlap. With equal priors, the minimum average error is one-half of the overlap, as in Eq.~\eqref{eq:bayes-risk}. The exact log-likelihood ratio retains the information needed for the optimal binary decision; a cruder score may lose separation but cannot create it.}
\label{fig:tv}
\end{figure}

We can visualize this idea by projecting the enormous discrete space \(\mathcal X\) onto a one-dimensional evidence score. The natural ideal score is the log-likelihood ratio
\[
s(x)=\log\frac{P_A(x)}{P_H(x)}.
\]
Here we use the usual conventions that the score is \(+\infty\) when only the AI source assigns positive probability to \(x\), and \(-\infty\) when only the human source does. Positive finite values favor the AI source; negative finite values favor the human source. Under our equal-prior assumption, the optimal decision boundary lies at \(s(x)=0\). Figure~\ref{fig:tv} shows the resulting overlap schematically.

Better algorithms still matter. More training data can bring an implemented detector closer to the optimal decision rule, and a better 
model may uncover information that a crude detector misses. But there is a ceiling. No algorithm can beat the error imposed by the
 statistical experiment itself.

Equation~\eqref{eq:bayes-risk} does not say that detection must always fail. It describes the information available from a particular
observation. Changing the observation therefore changes the experiment. A longer text, for example, may contain more evidence. If 
the longer observation includes everything contained in the shorter one, giving the detector the additional text cannot make the best
 possible classifier worse. The separation between the two sources, however, need not increase strictly, and 
Eq.~\eqref{eq:bayes-risk} alone does not guarantee that the error will approach zero as more text is observed. That requires additional
assumptions about how human and machine texts differ as more evidence is collected.

Recent theory illustrates these two sides of the problem. Sadasivan and colleagues~\cite{Sadasivan2025} connect optimal discrimination performance 
to total variation and show experimentally that paraphrasing can substantially reduce detectability. Chakraborty and colleagues~\cite{Chakraborty2024} 
study the complementary question: when a distributional distinction persists, how much additional evidence---through multiple samples or longer 
sequences---is required to make reliable discrimination possible?

There is no contradiction. One result emphasizes what happens when the two distributions are pushed closer together; the other studies what additional 
evidence can extract from a separation that remains. Detection therefore succeeds or fails relative to both the distributions being compared and the amount 
and structure of the evidence supplied to the detector. In practice, neither stays fixed for long.

\section*{There is no canonical distribution of human writing}
The two-lottery picture is useful, but it hides much of the complexity of actual writing. There is no canonical, context-independent distribution \(P_H\)
of human prose.

A mathematician proving a lemma, a schoolchild writing under examination, a novelist revising the fifteenth version of a paragraph, and a scientist 
writing carefully in a second language are not plausibly samples from one stable statistical source. Add lawyers, journalists, doctors, businesspeople, 
theologians, and editors, and the simplification becomes harder to defend. Genre matters, as do subject, age, education, native 
language, editorial conventions, and individual writing habits.

Nor is there a unique \(P_A\). Different AI models produce different distributions, and even the output of a particular model depends on its version, 
prompt, system instructions, sampling parameters, fine-tuning, retrieval context, translation, and subsequent revision.

A more realistic formulation therefore replaces two fixed distributions with two families:
\[
P_H\in\mathcal P_H
\qquad\text{versus}\qquad
P_A\in\mathcal P_A.
\]
This changes the problem in an important way. It is no longer enough for one particular human distribution to be well separated from one particular 
AI distribution. A detector intended for general use must remain reliable across many different pairs drawn from these two families.

This is where the Liang study becomes instructive. Liang and colleagues traced much of the disparity they observed to \emph{perplexity}, a signal used directly or indirectly by several early AI-text detectors. 
Roughly speaking, perplexity measures how surprising a sequence of words appears to a particular language model. If the model assigns relatively high probability 
to the words that actually occur, the text has low perplexity; if the sequence repeatedly surprises the model, its perplexity is higher.

If machine-generated prose tends to have lower perplexity than the human writing used for comparison, then low perplexity can become useful 
evidence for machine generation. But humans can also write predictable prose. Someone writing in a second language may rely on a narrower 
vocabulary or on familiar grammatical constructions. Technical disciplines may encourage standardized formulations. A careful writer may suppress 
precisely the lexical eccentricities that would make a language model more surprised. The statistical signal may therefore be real without its cause being unique.

Liang and colleagues found that the TOEFL essays classified as AI-generated by all seven detectors had significantly lower perplexity. They then asked 
ChatGPT to enrich the vocabulary of the same essays. Perplexity increased, while the reported average false-positive rate fell sharply \cite{Liang2023}.
The intervention therefore increased actual AI involvement while moving the essays away from a feature on which the detectors appeared to rely.

This does not mean that the detectors had no statistical basis. On the contrary, at least part of their behavior tracked a real statistical property of the 
text: predictability under a language model. The problem is that this property was not uniquely diagnostic of how the text had been produced.

There is, however, an important counterpoint. Jiang and colleagues~\cite{Jiang2024} assembled a large collection of Graduate Record Examination 
essays and developed detectors specifically for that setting using linguistic and perplexity-based features. They reported near-perfect 
discrimination between their human-written and ChatGPT-generated samples and found no evidence that their detectors disadvantaged non-native 
English writers.

The two studies need not contradict one another. They describe different statistical experiments. Liang and colleagues evaluated widely available 
detectors on human essays drawn from particular existing collections. Jiang and colleagues constructed and evaluated detectors within a comparatively 
standardized assessment domain using data sampled for that purpose.

The contrast illustrates a general problem of statistical generalization. A classifier can perform extremely well on the population and conditions 
for which it was developed yet deteriorate when the population, genre, model, prompt, or writing process changes. Validation results therefore 
belong not merely to an algorithm, but to an algorithm evaluated under a particular data-generating process.

For a fixed pair \(P_H\) and \(P_A\), the previous section asked how much the two distributions overlap. Once we admit families \(\mathcal P_H\) and 
\(\mathcal P_A\), a harder question appears: how close can some plausible human distribution come to some plausible AI distribution? Excellent 
separation for one pair tells us little about the pairs that were never tested.

AI detection makes this problem especially difficult because the distributions themselves evolve. Models are updated or replaced. Prompting 
practices change. Human writers adapt to AI tools. Machine-generated passages are edited by people and may then be edited again by another 
model. What counted as the statistical distinction between ``human'' and ``AI'' in one setting may not be the same distinction a detector encounters 
later. The detector does not face two fixed targets. It faces two moving families.

\section*{What does 98\% mean?}
Let us return to the AI-detector number on the instructor's screen:
\[
98\%\ \text{AI-generated}.
\]
Most readers will naturally interpret this as a posterior probability,
\[
\Pr(A\mid X=x)=0.98:
\]
given this particular essay, there is a 98\% probability that it came from an AI source.

However, the percentage does not necessarily represent a posterior probability. It might be a transformed model score, a percentile relative 
to some reference collection, an estimated fraction of sentences judged AI-like, a similarity measure, or a quantity calibrated on a population 
quite different from the one now being tested. The percentage sign alone does not tell us which statistical object is being reported.

Even if the human and machine likelihoods were known exactly, a posterior probability would still require one more ingredient: prior probabilities \(\pi_A\) and \(\pi_H\) for the two sources, assigned before the essay is examined. Bayes' theorem gives
\begin{equation}
\Pr(A\mid x)
=\frac{\pi_A P_A(x)}{\pi_A P_A(x)+\pi_H P_H(x)}.
\label{eq:bayes}
\end{equation}
In the population-based example that follows, these prior probabilities encode the \emph{base rate}: how common AI-generated essays are in the population under consideration before we
inspect the particular essay. This matters because even a detector that distinguishes the two groups quite well can produce many false alarms if
genuine AI-generated essays are rare.

Before returning to a numerical score such as 98\%, consider the simpler case in which an imaginary detector reports only a flag: positive or negative. Suppose that when an essay really is AI-generated, the detector flags it 95\% of the time. This quantity,
\[
\Pr(+\mid A)=0.95,
\]
is its \emph{sensitivity}. Suppose also that when an essay really is human-written, the detector correctly leaves it unflagged 95\% of the time. 
This quantity,
\[
\Pr(-\mid H)=0.95,
\]
is its \emph{specificity}. Equivalently, its false-positive rate is 5\%:
\[
\Pr(+\mid H)=0.05.
\]

These numbers may sound excellent. But they do not yet answer the instructor's question. They tell us how the detector behaves when the source 
is already known. The instructor faces the reverse problem: after seeing a positive result, how likely is it that the essay actually came from AI\@?

Suppose the detector is used in a population where only 1\% of essays are AI-generated. Imagine 10,000 essays from that population. About 
100 would be AI-generated. With 95\% sensitivity, approximately 95 of those would be flagged. The remaining 9,900 essays would be human-written, and a 
false-positive rate of 5\% would incorrectly flag about 495 of them. The detector would therefore produce approximately
\[
95+495=590
\]
positive results, of which only 95 would correspond to AI-generated essays. Thus
\[
\Pr(A\mid +)
\approx
\frac{95}{95+495}
\approx 0.16.
\]
Under these illustrative assumptions, the detector's positive predictive value---the probability of AI generation after a positive result---would be only about 16\%.

Nothing about the detector has changed. Its sensitivity is still 95\%, and its specificity is still 95\%. What changed is the frequency of 
AI-generated essays in the population in which the detector is being used.

To see how strongly this matters, suppose instead that half of all essays are AI-generated. Among 10,000 essays, there would then be 
5,000 AI-generated and 5,000 human-written. The same detector would flag approximately 4,750 of the AI essays and falsely flag 250
 of the human essays. Now
\[
\Pr(A\mid +)=\frac{4750}{4750+250}=0.95.
\]
The positive predictive value of the same detector has gone from about 16\% to 95\% simply because the population base rate has changed.

This is the base-rate problem in its simplest form. Sensitivity tells us
\[
\Pr(+\mid A),
\]
the probability of a positive result assuming that the essay came from AI\@. What the instructor wants to know is
\[
\Pr(A\mid +),
\]
the probability that the essay came from AI given that the detector returned a positive result. These conditional probabilities point in 
opposite directions. Bayes' theorem connects them, but it does not make them equal. Two decimal places do not make them the same quantity.

This also exposes a practical difficulty with a number such as \(98\%\). To interpret it as a posterior probability, we would need to know 
not only how the detector responds to human and AI text, but also what prevalence of AI-generated text was assumed or represented during its calibration. 
A detector calibrated in a collection containing equal numbers of human and AI texts need not produce calibrated posterior 
probabilities in a classroom where AI-generated submissions are much rarer---or much more common.

Bayes' theorem also shows what certainty about a particular essay would demand. If both prior probabilities are positive, then
\[
\Pr(A\mid x)=1
\]
requires
\[
P_H(x)=0
\qquad\text{and}\qquad
P_A(x)>0.
\]
Within this simplified two-source model, the observed essay would have to be impossible under the modeled human distribution. It is not enough for the essay to be 
unusual for a human writer, unusually polished, low in perplexity, or even much more probable under the AI model. As long as
\(P_H(x)>0\), the posterior remains strictly below one. If both \(P_H(x)\) and \(P_A(x)\) are zero, the model has assigned zero probability to the observation and the posterior is undefined; that is a failure of the model, not certainty about the source. For ordinary prose, assigning probability zero to human production is an audacious claim.

A likelihood ratio much greater or much smaller than one can nevertheless constitute powerful evidence. It can move a reasonable prior 
probability far toward one source or the other. But evidence is not certainty, and a score displayed as a percentage does not acquire a 
probabilistic interpretation merely because it looks like one.

\section*{Turing could ask another question}
\begin{figure}[H]
\centering
\includegraphics[width=0.96\textwidth]{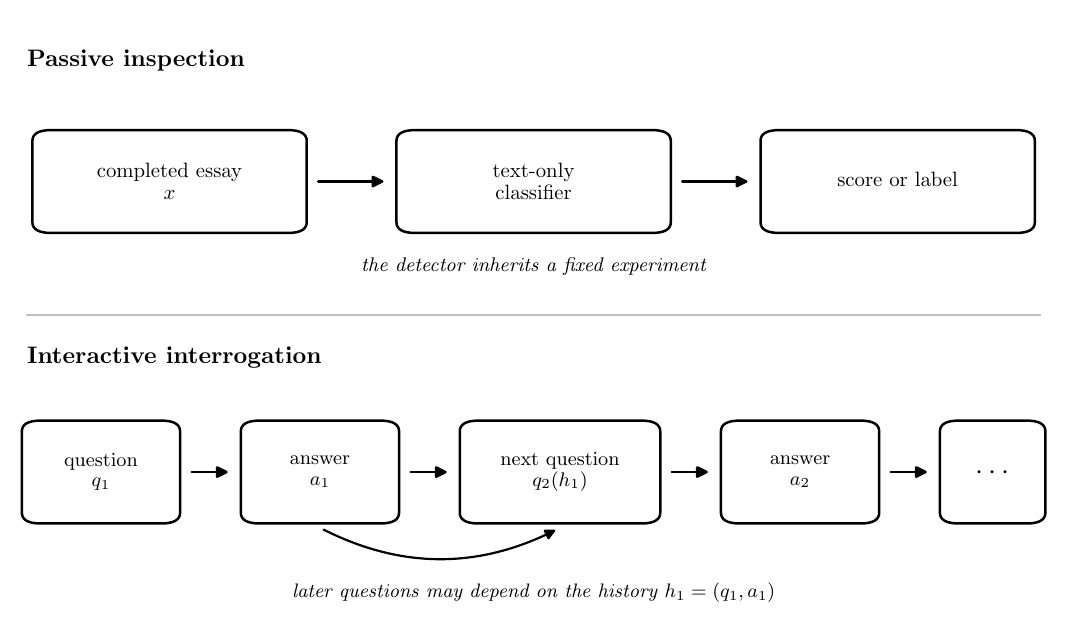}
\caption{\textbf{Inspection and interrogation.} In passive detection, a completed essay is passed to a classifier, which returns a 
label or score. In an interactive test, the interrogator chooses a question, observes an answer, and chooses the next question 
in light of the exchange so far. The detector inherits an experiment; the interrogator helps design one.}
\label{fig:interrogation}
\end{figure}

The preceding sections have treated the essay as a fixed object. Once it has been submitted, the detector receives whatever
information the text happens to contain and must work with that evidence alone. Socrates made the same complaint about written 
words: ask them a question and they say the same thing again \cite{Plato1997}. Turing's imitation game suggests another possibility: 
change the experiment. Figure~\ref{fig:interrogation} contrasts passive inspection with an adaptive interrogation.

In the imitation game, the interrogator is not restricted to inspecting a completed transcript. She can ask a question, observe 
the answer, and decide what to ask next. A surprising answer can be pursued. An inconsistency can be revisited. A line of 
questioning that proves uninformative can be abandoned.

We can express this in the same probabilistic language used earlier. Fix the number of rounds and the set of allowable questions. 
Let \(\pi\) denote a questioning strategy: a rule for choosing each new question from the exchange observed so far. Once the human 
and machine response mechanisms are fixed, a strategy \(\pi\) induces one distribution over complete transcripts when the respondent 
is human and another when the respondent is a machine:
\[
P_H^\pi
\qquad\text{and}\qquad
P_A^\pi.
\]
For any fixed strategy, the total-variation identity from Equation~\eqref{eq:bayes-risk} still applies. With equal prior probabilities, the 
smallest possible classification error after that interrogation is
\[
R_{\mathrm{best}}(\pi)=\frac12
\left[
1-\operatorname{TV}\!\left(P_H^\pi,P_A^\pi\right)
\right].
\]

But something has changed. The interrogator now has some control over the distributions themselves. Instead of accepting 
whatever evidence happens to be present in a finished essay, she can choose questions intended to make the two possible
sources behave differently. In this abstraction, the best questioning strategy seeks
\[
\sup_{\pi}
\operatorname{TV}\!\left(P_H^\pi,P_A^\pi\right).
\]

Equivalently, the smallest achievable error becomes
\[
\inf_{\pi}R_{\mathrm{best}}(\pi)=\frac12
\left[
1-
\sup_{\pi}
\operatorname{TV}\!\left(P_H^\pi,P_A^\pi\right)
\right].
\]
This is not Turing's notation. It is experimental design written in probabilistic form.

Adaptivity can matter because the next observation need not be chosen in advance. A fixed questionnaire asks the same 
sequence regardless of the answers. An adaptive interrogator can make question \(k+1\) depend on what happened in the 
first \(k\) exchanges. Within these constraints, nonadaptive strategies are available as special cases.. Allowing adaptive
questioning therefore cannot make the best possible experiment less informative; in favorable circumstances, it can make 
it substantially more so.

What might this mean in the classroom? An instructor could ask why a particular example was chosen, request that an 
argument be restated in different terms, introduce an unexpected assumption, return to an earlier claim, ask where a quotation 
or factual assertion came from, or require the reasoning to survive a new example. An answer that raises a new question can 
itself determine what is asked next. A detector examining a completed essay cannot do this; its evidence was 
fixed before the analysis began.

The distinction is important because the instructor may not ultimately care about the same question as the detector. 
A detector asks whether a static text resembles one statistical source more than another. An instructor may instead want 
to know whether a particular student understands the argument, can reproduce its reasoning, or can extend it under new 
conditions.

An oral examination, a supervised piece of writing, or a request to modify and defend an argument does not perfectly 
reconstruct the history of the submitted document. Nor is such evidence infallible: performance can depend on time 
pressure, verbal fluency, familiarity with the format, and many other factors. But these procedures generate new
observations of the student under controlled conditions rather than attempting to infer an entire history from the 
prose alone. That is a different statistical experiment. Sometimes a weak experiment needs a better experiment, not a more confident classifier.

\section*{What exactly is AI-written?}
There is a more fundamental difficulty with AI detection, one that arises before questions of
accuracy. We may not have clearly defined the property that the detector is supposed to identify.

Consider several writing processes. One writer drafts an essay independently and asks an AI system
only to correct the grammar. Another accepts an AI-generated outline but writes every sentence. A
third supplies the argument, receives a machine-generated draft, and then reconstructs every
paragraph while correcting the reasoning. A fourth develops the essay through an extended dialogue
with an AI system, accepting some suggestions, rejecting others, and allowing each response to shape
the next prompt. In another case, the model supplies ideas while the person supplies the language;
reverse those roles and the compositional history changes again.

Which of these finished documents should count as \emph{AI-written}? We often call such examples ``boundary
cases,'' but that expression already assumes that there is a single boundary to be found. There may
instead be several different dimensions: who supplied the ideas, who determined the structure, who
produced the words, who revised them, and how dependent the final work was on the assistance
received. Different policies may reasonably place the boundary in different places.

Measurement science gives us a useful analogy. Before evaluating a measurement procedure, one
specifies the \emph{measurand}: the quantity intended to be measured. Assigning a text to the class
``AI-generated'' concerns a nominal property, not a quantity; in the terminology of the International
Vocabulary of Metrology, classification of such a property is not itself measurement \cite{JCGM2012}.

A detector may nevertheless calculate a numerical score and use it to assign a class. The score
needs an interpretation; the class needs an operational definition. Reporting a percentage does not
by itself establish that a quantity has been measured, nor does it define the writing process that
the label denotes. Only after deciding what counts as ``AI-generated'' can quantities such as
sensitivity, specificity, and error rates be interpreted.

But several quite different definitions are possible. Does ``AI-generated'' mean that any
machine-produced words remain in the final document? That most of the words came from a model? That
the first draft was generated by AI\@? That the organization or principal argument originated with
the machine? Or does it mean that the writer could not have produced the work without AI assistance?
These definitions concern different properties of the writing process, and there is no reason for
them to agree. If ``AI-written'' is intended to describe how a document was produced, then
provenance---the history of its production---is not a separate topic. The intended target is some
feature of that provenance.

A document may, for example, have a compositional history such as
\[
\text{human}
\longrightarrow \text{AI}
\longrightarrow \text{human}
\longrightarrow \text{AI}
\longrightarrow \text{final text}.
\]
The final classification compresses this entire history into a single binary label---``human'' or
``AI.'' Assigning the label does not resolve the ambiguity in the process; it merely hides that
ambiguity inside the definition of the class.

We can make the point slightly more formally. Let \(Z\) denote the complete, unobserved history by
which a document was produced, and let \(X\) denote the finished text. A dataset labeled ``AI'' and
``human'' implicitly applies some rule
\[
Y=g(Z)
\]
that turns the complicated history \(Z\) into a binary label \(Y\). A classifier then attempts to
infer \(Y\) from \(X\). Its accuracy can therefore be excellent even when the labeling rule \(g\)
has little to do with the question its eventual user wants answered. Before asking for sensitivity,
specificity, or accuracy, we should therefore ask a preliminary question: what exactly counted as
``AI'' in the training and validation data, and does that definition correspond to the decision now
being made?

Several distinct inferential tasks are easily conflated here. Turing's question concerns
\emph{distinguishability}: can observable machine behavior be distinguished from human behavior? An
AI-text detector is ordinarily used for \emph{attribution}: given a completed text, which proposed
source or production process better explains the observation? A university, publisher, court, or
reader may instead require \emph{authentication}: what evidence connects this particular work to
this particular person? In many practical settings there is also a question of \emph{compliance}:
was the work produced under the conditions that were permitted? A student may unquestionably be the
author of an essay and still have used prohibited assistance; conversely, a student may have used an
allowed grammar tool without surrendering authorship of the work.

These distinctions matter because a text's source and production history are not directly observable
in the finished text. We see the words, not the history that produced them. To move from ``this
passage has features associated with machine-generated text'' to ``this student used AI,'' and then
to ``this student violated the rules,'' is to make a succession of increasingly strong inferences:
\[
\begin{aligned}
\text{statistical resemblance}
&\longrightarrow \text{class-level source attribution}\\
&\longrightarrow \text{a claim about a particular person's actions}\\
&\longrightarrow \text{a claim of noncompliance}.
\end{aligned}
\]
The arrows are not logical implications. Each marks a separate inferential step requiring additional
evidence or assumptions. Attribution is therefore an evidentiary commitment: a conclusion about an unobserved history, not a
property read directly from prose.

An AI detector may contribute evidence at the first one or two stages, depending on what it was
trained and validated to distinguish. It cannot, from stylistic resemblance alone, establish a
document's complete history: human prose may resemble model output, and heavily AI-assisted prose
may be revised until it resembles ordinary human writing.

Watermarks and cryptographically signed provenance records illustrate two ways to change the
experiment rather than escape Equation~\eqref{eq:bayes-risk}. A watermark modifies \(P_A\) by embedding a
detectable signal in machine-generated text; a signed provenance record augments the observation,
so that the verifier receives evidence about the text's production alongside the text itself \cite{Kirchenbauer2023,C2PA2026}.
A valid signature does not by itself establish the truth or completeness of the recorded history.
Neither method reconstructs a complete production history from unauthenticated prose alone.

This suggests a different strategy when provenance genuinely matters. Rather than trying to
reconstruct an entire writing process after the fact from stylistic traces, we can preserve evidence
while the work is being produced. Draft histories, contemporaneous notes, source records,
authenticated workspaces, supervised writing, and oral explanation are all imperfect, and some can
be manipulated. Some of these approaches preserve evidence of the production process; while others 
provide additional evidence of the writer’s understanding and capabilities. They can therefore supply 
information unavailable from the finished prose alone.

\section*{Who wrote this?}
The paper you are reading makes the problem unusually literal. Its own provenance is mixed.

I am a mathematical scientist whose work involves measurement in the chemical and life sciences. Questions 
about what an instrument, algorithm, or mathematical procedure actually measures—and under what 
assumptions its output can be interpreted—are therefore difficult for me to leave alone. After examining 
the outputs of AI-text detectors, I realized that it was far from clear what these systems were actually 
measuring.

This paper began with that question. In September 2026, I asked OpenAI's GPT-5.6 Pro whether an AI-written 
essay could ever be identified with certainty, and I suggested connecting the problem to Turing's imitation 
game and to the ``reverse'' Turing test posed by AI-text detection. GPT-5.6 Pro connected the question to 
binary hypothesis testing and statistical decision theory, proposed early structures and titles, identified 
candidate literature, and generated portions of early drafts.

I chose the problem, journal, audience, and eventual direction of the argument. I decided which examples 
and results should remain central, including the decision to retain the Liang study while confronting it with 
important counterevidence. I checked the cited literature, rejected or corrected claims that I judged inaccurate, 
corrected mathematical and conceptual errors, and reorganized the argument.

In revising the manuscript, I rewrote most sections and added substantial explanations. I occasionally asked 
GPT-5.6 Sol to provide language-editing assistance for sections I had rewritten. I created the figures and 
made the final decisions about the mathematics, interpretation, sources, organization, and wording. I take 
full responsibility for the resulting manuscript.

Before submitting the manuscript for publication, I also asked GPT-5.6 Pro and Anthropic Claude Fable 5.1 
to review it critically and identify issues that reviewers might raise. I evaluated the issues they identified and 
made corrections where I judged them warranted.

Calling it ``coauthored'' would introduce a different confusion. Mathematical tools have long made
substantive contributions without thereby becoming authors. \texttt{Maple}\texttrademark{} or another 
computer algebra system, for example, may carry out a complex and interactive symbolic calculation 
on which a perturbation analysis depends, while the mathematician chooses the problem, determines the 
calculation's relevance, interprets its output, and accepts responsibility for the resulting claim. The point is not
that the software contributed nothing. It is that contribution and authorship are different roles.

Generative AI can contribute more broadly---through language, structure, criticism, examples, or
proposed ideas---but breadth of contribution does not by itself establish authorship. Scholarly
authorship also carries accountability for claims, sources, errors, judgments, and the decision to
publish; that accountability remains mine.
A recent disciplinary statement articulates a closely related position.
The Leiden Declaration on Artificial Intelligence and Mathematics calls
for transparent disclosure of automated-tool use, assigns responsibility
for correctness and citation to human authors, and does not treat automated
systems as authors \cite{LeidenDeclaration2026}.

A provenance record therefore says something that a binary label cannot. The history of this
manuscript is not
\[
\text{human}
\qquad\text{or}\qquad
\text{AI}.
\]
It is closer to
\[
\begin{aligned}
\text{human question}
&\longrightarrow \text{AI response}
\longrightarrow \text{human judgment}\\
&\longrightarrow \text{AI criticism}
\longrightarrow \text{human revision}\\
&\longrightarrow \cdots
\longrightarrow \text{final manuscript}.
\end{aligned}
\]

That history is precisely the kind of object a text-only detector is being asked to reconstruct from
the final prose. It may detect statistical traces associated with one stage of the process. It may
even provide useful evidence about some possible histories. But the finished sequence of words does
not contain a complete transcript of how the document was conceived, argued over, corrected, and
made. The example also shows that mixed provenance need not imply mixed accountability:
contributions can be distributed across people and tools, while responsibility for this manuscript
remains mine.

So who wrote this paper? The answer is available because its history has been recorded rather than
inferred from its style.

The title can therefore be answered, but not with one bit.

\clearpage
\backmatter

\bmhead{Acknowledgements}The author used OpenAI's ChatGPT (GPT-5.6 Pro and GPT-5.6 Sol; accessed 
September 2026) and Anthropic Claude (Fable 5.1; accessed September 2026) as described above. The author 
critically evaluated, revised, and verified material produced with their assistance and takes full responsibility 
for the accuracy, interpretation, and integrity of the final manuscript.  Both figures were generated 
programmatically and were not produced using a generative-image system.

\bmhead{Declarations}\par

 \emph{Competing interests:} The author has no relevant financial or non-financial interests to disclose. 
\emph{Data and code availability:} No empirical data were generated or analyzed.

\bibliography{references}

\end{document}